\documentclass[lettersize,journal]{IEEEtran}
\usepackage{amsmath,amsfonts}
\usepackage{algorithm}
\usepackage{array}
\usepackage[caption=false,font=normalsize,labelfont=sf,textfont=sf]{subfig}
\usepackage{textcomp}
\usepackage{stfloats}
\usepackage{url}
\usepackage{color}
\usepackage{verbatim}
\usepackage{graphicx}
\usepackage{cite}
\usepackage[margin=2.5cm]{geometry} 
\usepackage{graphicx} 
\usepackage{multirow}
\usepackage{lscape}
\usepackage{ntheorem}
\usepackage{hologo}
\usepackage{amssymb}
\usepackage{calc}
\usepackage{ntheorem}
\usepackage{bm}
\usepackage{hyperref}
\usepackage{mathrsfs}
\usepackage{algpseudocode}
\usepackage{booktabs}
\allowdisplaybreaks[4]
\newtheorem{remark}{Remark}
\newtheorem{assumption}{Assumption}
\newtheorem{lemma}{ Lemma}
\newtheorem{theorem}{ Theorem}
\newtheorem{corollary}{ Corollary}
\begin{document}

\title{Judicial Sentencing Prediction 
Based on Hybrid Models and Two-Stage Learning Algorithms}

\author{ Ruifen Dai, Xin Zheng, Fang Wang, and Lei Guo, \IEEEmembership{Fellow, IEEE}
\thanks{This paper was supported by the National Natural Science Foundation of China under Grant Nos.  T2293773,
72371145 and 12288201, the National Key Research and Development Program under Grant No. 2024YFC3307200, the Special Funds for Taishan Scholars Project of Shandong Province, China, under
 Grant No. tsqn202211004. 
 (Corresponding Author: Fang Wang, Lei Guo.)
 }
\thanks{
Ruifen Dai and Fang Wang are with  the Data Science Institute, Shandong University, Jinan 250100, China. (e-mails: dairuifen@mail.sdu.edu.cn, wangfang226@sdu.edu.cn).

Xin Zheng and Lei Guo are with the State Key Laboratory of Mathematical Sciences, Academy of Mathematics and Systems Science, Chinese Academy of Sciences, Beijing 100190, China, and also with the School of Mathematical Sciences, University of Chinese Academy of Sciences, Beijing 100049, China. (e-mails: zhengxin2021@amss.ac.cn, lguo@amss.ac.cn).}
}

%
%
\maketitle

\begin{abstract}
The investigation of legal judgment prediction (LJP), such as sentencing prediction, has attracted broad attention for its potential to promote judicial fairness, making the accuracy and reliability of its computation result an increasingly critical concern. In view of this, we present a new sentencing model that shares both legal logic interpretability
and strong prediction capability by  introducing a two-stage learning algorithm. Specifically, we first construct a hybrid  model that synthesizes a mechanism model based on the main factors for sentencing with a neural network modeling possible uncertain features.
We then propose a two-stage learning algorithm: First, an  adaptive  stochastic gradient (ASG)   algorithm is used to get good estimates for the unknown parameters  in the mechanistic component of the hybrid model. Then, the Adam optimizer tunes all parameters to enhance the predictive performance of the entire hybrid model. The asymptotic convergence of the ASG-based adaptive predictor is established without requiring any excitation data conditions, thereby providing a good initial parameter estimate for prediction.  Based on this,  the fast-converging Adam optimizer further refines the parameters to enhance overall prediction accuracy. 
Experiments on a real-world dataset of intentional injury cases in China show that our new hybrid model combined with our two-stage ASG-Adam algorithm, outperforms the existing related methods in sentencing prediction performance, including those based on neural networks and saturated mechanism models.
\end{abstract}


\section{Introduction}
\IEEEPARstart{I}{n} recent years, the investigation of various problems in legal judgment prediction (LJP), such as sentencing prediction, has attracted broad research attention due to its vital role in enhancing fairness in the judicial system~(\cite{1}-\cite{40}). Among these problems, ensuring the high reliability and accuracy of sentencing prediction results is an important requirement in judicial practice. 
This requirement calls for sentencing models that integrate the sentencing logic for interpretability and capture case-specific uncertainties not specified explicitly  in law for prediction capability.
However, such models have been rarely developed in the literature. 
Moreover, most existing theories on the current mainstream gradient-based learning algorithms are not directly applicable to sentencing prediction, since they rely on common yet stringent statistical assumptions on the data, such as the independent and identically distributed (i.i.d.) data assumption, which are not the case for real-world judicial datasets. 
Therefore, developing sentencing prediction models based on judicial logic and establishing the performance guarantees  of the  associated learning algorithms, hold significant value in both theory and practice.




Although most existing LJP studies have not adequately addressed the aforementioned challenges, their efforts toward achieving high-precision sentencing prediction offer valuable insights (see, e.g., \cite{41}–\cite{35}, \cite{11}–\cite{42}, \cite{12}, \cite{21}), which can be broadly summarized in terms of model construction and algorithm design in the following:

First, extensive research has been conducted on using machine learning models for sentencing prediction. For example, \cite{41} introduced a convolutional neural network (CNN)-based model for sentencing classification prediction. To address the need to model subtask dependencies, \cite{1} employed a multi-task learning framework  TopJudge that represents the dependencies among LJP subtasks using a directed acyclic graph. To further capture the relationships between subtasks, \cite{add31} proposed a multi-perspective bi-feedback network model enhanced by a word collocation attention mechanism, which reflects subtask dependencies through a topology-aware design. Subsequently, \cite{32} modeled LJP as a node classification problem over a global consistency graph. \cite{33} presented the LADAN model, which combines a graph neural network (GNN) and an attention mechanism to distinguish confusing law articles for more accurate LJP. 
In addition, \cite{34} introduced a NeurJudge model that sets BERT (see \cite{45}) as a judgment encoder and leverages intermediate subtask results to partition fact descriptions into circumstances for guiding subsequent subtask predictions.
\cite{35} presented a QAjudge model composed of a question net and an answer net to visualize the prediction processes. 

Despite the widespread application of the machine learning sentencing models, they fail to fully incorporate the sentencing logic in the sentencing guidelines.
In view of this, \cite{11} proposed an interpretable saturated mechanism model based on the sentencing logic of ``starting point—benchmark sentence—announced sentence", and applied it to sentencing prediction for the crime of intentional injury. To improve sentencing prediction accuracy, \cite{42} incorporated saturated boundaries to ensure that neural network–predicted sentences fall within the statutory sentencing range. \cite{Guo2024} and \cite{newmodel} refined the sentencing mechanism model proposed in \cite{11} by incorporating the temporal logic underlying sentencing adjustments for benchmark sentences. The  impact that not explicitly encoded in sentencing guidelines, although potentially varying across different cases, is coarsely represented by a fixed constant in \cite{Guo2024} and \cite{newmodel}.

Second, many efforts have been devoted to designing algorithms for sentencing prediction. For instance,
\cite{12} proposed a two-step Newton-type adaptive algorithm and analyzed the related performances for both  parameter estimation and sentencing prediction, without resorting to the traditional persistence of excitation (PE) on the system data.
Subsequently, \cite{21} designed a more robust two-step weighted $l_1$-based Newton-type algorithm to improve prediction performance and established its global convergence.
Both algorithms in \cite{12}-\cite{21} were applied to sentencing prediction tasks based on the saturated mechanism model proposed in \cite{11}. However,  all the above algorithms incur high computational costs when processing large-scale and high-dimensional judicial data due to their second-order nature.
Different from  the Newton-type  algorithms, the gradient descent algorithms and their variants, such as Adam \cite{43}, Adadelta \cite{44}, and AdamW \cite{add46}, are also used for sentencing prediction  based on neural network models (e.g., \cite{41}-\cite{35}, \cite{42}).  However, 
neural networks are just mathematical approximation and  can hardly provide judicial interpretability for the sentencing results.

To overcome the aforementioned shortcomings, we establish a new sentencing
model that possesses both high legal interpretability and strong predictive capability by introducing a two-stage gradient algorithm. The main contributions can be summarized as follows:
\begin{itemize}

    \item We will introduce a hybrid model that integrates the sentencing-logic-based mechanism model with a neural network, which can not only reflect the main factors for sentencing, but also can reflect the influence of possible uncertain factors not explicitly encoded in the sentencing guidelines.

    \item  We will propose a two-stage learning (TSL)  algorithm for the prediction of the hybrid sentencing model. To be specific, an adaptive stochastic gradient (ASG)  algorithm will provide a good initial value of the mechanistic component parameters for the Adam optimizer in the next stage to 
fine-tune all parameters to improve the sentencing prediction performance.
    
    \item We will rigorously establish the global asymptotic convergence of the ASG-based adaptive predictor without requiring any excitation data conditions. This theoretical guarantee ensures that the ASG-based initialization provides a reliable approximation to globally optimal prediction, thereby offering a solid basis for the refined prediction by the Adam optimizer.
    
     \item Empirical experiments using a real-world judgment dataset for the crime of intentional injury (CII) will show that our proposed hybrid model and two-stage algorithm can achieve  high accuracy in sentencing prediction, outperforming other known related methods.
\end{itemize}

The structure of this paper is as follows: Section II will introduce the new sentencing model. Section III will present the two-stage ASG-Adam algorithm. Section IV will provide the global asymptotic convergence theory of the ASG algorithm whose proof will be shown in the appendix. In Section V, we will demonstrate the advantages of the proposed model and algorithm using a real-world sentencing dataset. The concluding remarks will be presented in the final section.

\section{A Hybrid  SMNN  Model }\label{section2}

In this section, we propose a hybrid model that integrates the sentencing mechanism (SM) in Chinese Criminal Law with a neural network (NN) accounting for possible uncertainties in sentencing,  abbreviated  as SMNN model, as follows:
\begin{equation}\label{modelsimulation1}
	\begin{aligned}
		z_{k+1} = S_k \bigg(&\underbrace{\left[ a_k + b x_k^{(1)} + c x_k^{(2)} \right]}_{\text{benchmark sentence}} \\
		&\times \underbrace{\prod_{i=1}^{m_1}(1 + p_i  v_k^{(i)} )}_{\text{influence of  primary factors}}\\
		&\times \underbrace{\left[ 1 + \sum_{j=1}^{m_2} q_j u_k^{(j)} + e_k \right]}_{\text{influence of other factors}} + w_{k+1} \bigg),
	\end{aligned}
\end{equation}
where the time-varying bias term $e_k$
and the saturated function $S_k(\cdot)$ are  defined  respectively as follows:
\begin{equation}\label{simulation4}
	e_k = \Gamma \sigma\left( B \sigma\left( A \eta_k + b^{(1)} \right) + b^{(2)} \right) + b^{(3)},
\end{equation}
\begin{equation}\label{2}
	S_{k}(x) =
	\begin{cases}
		L_k, & x < L_k; \\
		x, &L_k \leq x \leq N_k; \\
		N_k,& x > N_k,
	\end{cases}
\end{equation}
 and where \( z_k \in \mathbb{R} \) is the pronounced sentence, \( a_k \in \mathbb{R} \) is the sentencing starting point, and \( x_k^{(i)} \in \mathbb{R}, i = 1, 2 \), represent the factors determining penalty amounts, e.g., in the CII \footnote{Related definition and interpretation for the crime of intentional injury can be found in Article 234 of Chinese Criminal Law, available on the website https://flk.npc.gov.cn/.}, they represent the number  of seriously and minorly injured victims, respectively. \( b, c \in \mathbb{R} \) quantify the additional sentence for the offender corresponding to a one-unit increase in \( x_k^{(i)} \). \(  v_k^{(i)}, i = 1, \cdots, m_1 \) denote primary sentencing factors with application priority, and \( u_k^{(j)}, j = 1, \cdots, m_2 \) denote other sentencing factors. \( p_i, i = 1, \cdots, m_1 \) and \( q_j, j = 1, \cdots, m_2 \) are unknown weighting parameters for \(  v_k^{(i)} \) and \( u_k^{(j)} \), respectively. $m_1$ and $m_2$ are the number of primary sentencing factors and other sentencing factors, respectively.
 \( w_k \in \mathbb{R} \) denotes potential random noise effects. 
  \( e_k \in \mathbb{R} \) in (\ref{simulation4}) is the bias term reflecting the comprehensive influence of possible factors not specified explicitly in the law, and 
 \( A \in \mathbb{R}^{m \times m_3}, B \in \mathbb{R}^{m \times m}, \Gamma \in \mathbb{R}^{1 \times m} \) are unknown weight parameter matrices or vectors, and \( b^{(1)} \in \mathbb{R}^{m}, b^{(2)} \in \mathbb{R}^{m}, b^{(3)} \in \mathbb{R} \) are unknown bias parameters in the neural network. Here, $m$ represents the number of neurons in each of the two hidden layers. \( m_3 \) denotes the number of   possible  factors that not specified in the  law. \(\eta_k \in \mathbb{R}^{m_3}\) denotes the vector of factors that not explicitly encoded  in the  law  of the \(k\)-th case's sentencing.
  \(\sigma(X)\) is an activation function, such as the Rectified Linear Unit (ReLU) function, defined as 
$\left( \max(x_1, 0), \max(x_2, 0), \dots, \max(x_m, 0) \right)^\tau$  with $x_i$
being the element of $X$ in the $i$-th row, $i = 1, \cdots , m$.
 \( L_k \in \mathbb{R} \) and \( N_k \in \mathbb{R} \) in (\ref{2}) are the lower and upper bounds for the announced sentence, respectively, which are prescribed in Chinese Criminal Law and may vary for different crimes and  criminal cases.

\begin{remark}
The SMNN model (\ref{modelsimulation1})  is mainly based on Chinese Criminal Law and sentencing guidelines\footnote{See https://flk.npc.gov.cn/ and https://www.court.gov.cn/.}, it specifically follows the sentencing logic of “starting point—benchmark sentence—announced sentence” and incorporates temporal adjustments for benchmark sentences (see \cite{Guo2024}, \cite{newmodel} for more details).
Moreover, the model employs a neural network to approximate the comprehensive influences of factors that are not explicitly encoded in the law. Unlike treating the bias term as a fixed constant in previous studies (\cite{Guo2024} and \cite{newmodel}), the neural network can approximate the case-specific influences, enhancing the model’s flexibility and predictive performance.

The SMNN model (\ref{modelsimulation1}) will degenerate to the mechanism model in \cite{11}, if we take $v_k^{(i)} =0, i=1,\cdots,m_1$ and replace $e_k$ with a constant $e$.

\end{remark}

\begin{remark}
The saturated function (\ref{2}) ensures that the sentences are restricted to the statutory penalty range, and when a sentence exceeds the upper limit \( N_ k \) or falls below the lower limit \( L_k \), the final judgment is capped at \( N_k \) or \( L_k \), respectively. 

Besides, the saturated phenomena described by (\ref{2}) are also commonly found in various application fields, including engineering systems (\cite{add15}-\cite{add16}), economic behavior analysis (\cite{add17}-\cite{add18}), biomedical systems (\cite{add19}), and others. These applications collectively emphasize the importance of saturated scenarios in complex nonlinear modeling and analysis.

\end{remark}


\begin{remark}
The SMNN model mainly comprises two parts: (i) a mechanistic component grounded in the sentencing logic, and (ii) a neural network that captures the influence of possible residual factors (RF) not explicitly encoded in sentencing law. 
Since the mechanistic component already encodes the factors and rules explicitly stated in the sentencing law, which plays a dominating role in sentencing, the influence of the RF may not be significant.
Therefore, we first assume the bias term $e_k=e$ to be a constant reflecting the averaged effect of the RF, in order to facilitate a more straightforward estimation of the mechanistic component. After this stage, the highly nonlinear neural network will be introduced to enhance the prediction accuracy for each case.
\end{remark}

To avoid non-convex optimization problems to get better
estimators of the mechanistic component parameters, we rewrite the internal structure of
(\ref{modelsimulation1}) as a linearly parameterized regression
by increasing the dimension of  both regression vectors and parameter vectors. Here the bias term $e_k =e$.
To be specific, set the expanded regressor and parameter vector as follows:
\begin{align}
	\phi_k = \big[ 
	& a_k\phi_{1k}^T,  a_k(\phi_{2k} \otimes \phi_{1k})^T, x_k^{(1)}\phi_{1k}^T, \notag \\
	& x_k^{(1)}(\phi_{2k} \otimes \phi_{1k})^T, x_k^{(2)}\phi_{1k}^T,  x_k^{(2)}(\phi_{2k} \otimes \phi_{1k})^T 
	\big]^T, \label{phi_k}\\
	\theta = \big[
	& (1+e)\vartheta_1^T, (\vartheta_2 \otimes \vartheta_1)^T, b(1+e)\vartheta_1^T, \notag \\
	& b(\vartheta_2 \otimes \vartheta_1)^T,  c(1+e)\vartheta_1^T, c(\vartheta_2 \otimes \vartheta_1)^T
	\big]^T, \label{simulation3}
\end{align}
where
\begin{align}
	\phi_{1k} =& [
	1, \ z_k^{(1)}, \ \cdots, \ z_k^{(m_1)}, \ z_k^{(1)}z_k^{(2)}, \ \cdots, \ \notag\\
	& z_k^{(m_1-1)} z_k^{(m_1)}, \ \cdots, \ z_k^{(1)} \cdots z_k^{(m_1)}]^T, \\
	\phi_{2k} =& [ u_k^{(1)}, \ u_k^{(2)}, \ \cdots, \ u_k^{(m_2)} ]^T,\\
	\vartheta_1 =& [
	1, \ p_1, \ \cdots, \ p_{m_1}, \ p_1 p_2, \ \cdots, \ p_{m_1-1} p_{m_1}, \ \notag\\
	& \cdots, \ p_1 \cdots p_{m_1} ]^T, \\
	\vartheta_2 =& [ q_1, \ \cdots, \ q_{m_2} ]^T,
\end{align}
and the Kronecker product of two vectors is defined as
$\mathbf{a} \otimes \mathbf{b} =
[a_1 \mathbf{b},  ~a_2 \mathbf{b}, ~\cdots  ,~a_m \mathbf{b}
].$
Then (\ref{modelsimulation1}) can be rewritten as follows:

\begin{equation}\label{simulation2}
	z_{k+1} = S_k(\phi_k^T \theta + w_{k+1}).
\end{equation}

For the above new parameterized model (\ref{simulation2}), the dimension of the unknown parameter \(\theta\) is \(2^{m_1}(3 + 3m_2)\), which grows exponentially with \(m_1\), resulting in a huge demand for computational resources to handle such high-dimensional inputs. Taking the CII example as detailed below, the dimension of unknown parameter vector is 565,248, rendering many Newton-type methods (e.g., \cite{12}, \cite{21}) computationally infeasible due to excessive memory and time complexity.

\section{A Gradient-Based Two-Stage Algorithm} \label{section3}
In this section, we propose the two-stage  gradient learning algorithm for the SMNN model (\ref{modelsimulation1}).
  
To better introduce the algorithm, we first introduce the following notations and assumptions.

\subsection{Notations and Assumptions}


\textbf{Notations.} Throughout the sequel, $\| \cdot \|_1$ and $\| \cdot \|$ denote the $1$-norm and Euclidean norm of vectors or matrices, respectively.
The maximum and minimum eigenvalues of a square matrix $X$ are denoted by $\lambda_ {\max}\{X\}$ and $\lambda_ {\min}\{X\}$, respectively.
Let \(\{ \mathcal{F}_k \}\) be a non-decreasing sequence of \(\sigma\)-algebras, along with the associated conditional expectation operator \(\mathbb{E}[\cdot \mid \mathcal{F}_k]\).
Moreover, for any two sequences $\{a_n\}, \{b_n\}$ with $b_n>0$,  $a_n=O(b_n)$ means that there exists a constant $N>0$ such that $|a_n|/b_n \leq N$ for all $n>0$, and $a_n =o(b_n)$ means that $a_n/b_n \rightarrow 0$ as $n \rightarrow \infty$.

\begin{assumption} \label{assumption1}
The regressor $\phi_k$ is $ \mathscr{F}_k$-measurable and bounded. Also, there is a known constant $L>0$ such that $ \|\theta\|_1 \leq  L$.
\end{assumption}
Under Assumption \ref{assumption1}, there exists a positive bounded sequence $\{M_k\}$  such that $\max \limits_{1 \leq i \leq 2^{m_1}(3 + 3m_2)} |\phi_{k,i}| \leq M_k$ for all $k \geq 0$, where $\phi_{k, i}$ is the $i$-th component of $\phi_k$ and $M_k$ is bigger than the natural constant $e$. 
The time‐varying bound $M_k$ will be used in the subsequent design of the algorithm.

\begin{assumption} \label{assumption2}
The thresholds $L_k$ and $N_k$ defined in (\ref{2}) are $\mathscr{F}_k$-measurable, uniformly bounded with respect to the sampling path, and strictly ordered, satisfying $L_k < N_k ~~\text{for all }k$.

\end{assumption}
Note that Assumption (\ref{assumption2}) is quite general for  ensuring the saturated property of the function $S_k(\cdot)$.

\begin{assumption} \label{assumption3}
	The random noise $\{w_k , \mathscr{F}_k\}$ is a martingale difference sequence, and satisfies
	\begin{align}\label{7}
		&\sup_{k \geq 0} E\big[|w_{k+1}|^{4} | \mathscr{F}_k\big] < \infty, \quad a.s.
	\end{align}
	And the conditional expectation function $G_k(x) \triangleq E[S_k(x+w_{k+1})|\mathscr{F}_k]$ is differentiable, and its derivation function $G_k^{'}(x)$ satisfies
	\begin{align}\label{81}
		& 0 <\underline{g} =\inf_{|x| < C, k \geq 0}G_k^{'}(x) \leq \sup_{|x|< C, k \geq 0}G_k^{'}(x) =\overline{g} < \infty,
	\end{align}
	where $C>0$ is any constant.
\end{assumption}
For convenience of analysis, we also introduce the following notations used in the subsequent sections:
\begin{align}
	&\tilde{\theta}_k = \theta - {\theta}_k,\label{24}\\
	&v_{k+1}  =  z_{k+1} - G_k(\phi_k^{\tau} \theta) , \label{26}\\
	&\psi_k = G_k(\phi_k^{\tau} \theta) - G_k\big(\phi_{k}^{\tau}  {\theta}_{k}\big), \label{27}  \\
	& \underline{g}_k   = \inf  \limits_ {|x|\leq \max\{M_kL, M_k\|\theta_k\|_1\}}G_k^{'}(x), \\
    & \bar{g}_k = \sup \limits_{ {|x|\leq \max\{M_kL, M_k\|\theta_k\|_1\}}}G_k^{'}(x),
\end{align}
where $\theta_k, k \geq 0$ is the estimate for $\theta$,  which is obtained by our algorithm introduced below.
\subsection{A Gradient-Based Two-Stage  Algorithm}
Now, we propose our two-stage gradient algorithm, the details of this algorithm are outlined in Algorithm \ref{SG-Adam} below.  

\begin{algorithm*}
	\caption{Two-Stage Learning  (TSL) Algorithm}
	\label{SG-Adam}
	\begin{algorithmic}[1]
        \State \textbf{Input}:  The  temporally ordered training dataset $D_1= \{\phi_k, z_{k+1}\}_{k=0}^{n-1}$ ,  the held-out  testing dataset $D_2= \{\phi_k, z_{k+1}\}_{k=1}^{n_2}$,       
        the arbitrarily chosen initial estimates $\theta_{0} \in \mathbb{R}^{p \times 1}$, $ \hat{\Gamma}_{0} \in \mathbb{R}^{1 \times m}$, $\hat{B}_{0} \in \mathbb{R}^{m \times m}$, $\hat{A}_{0} \in \mathbb{R}^{m \times m_3}$, $\hat{b}_{0}^{(1)} \in \mathbb{R}^m$, $\hat{b}_{0}^{(2)} \in \mathbb{R}^m$, 
        $\hat{b}_{0}^{(3)} \in \mathbb{R}$, the epoch count of $N$, the batch size  $\mathcal{T}$ ($\mathcal{T} \leq n$), hyper-parameters $ \alpha>1$, $\eta_1 >0$,  $\varepsilon >0$,  $\mu\in (0,1]$ , $\beta_1 \in (0,1)$  and $ \beta_2 \in (0,1)$. 
         \State \textbf{Output}:  The final parameter estimate $\hat{{\Theta}}_{\lfloor n/\mathcal{T} \rfloor}^{(N)} $, the predictive sentences for sentencing cases in $D_2$.
         
		\begin{footnotesize}
    \State \textbf{\#~Stage 1: ASG-based initialization~for~the~ unknown mechanistic parameter in the model (\ref{modelsimulation1})}.
    \State\hspace{2em} \textbf{for $k=0$ to $n-1$ do} 
		\begin{align}
    	\theta_{k+1} &= \theta_k  + \frac{ \mu \bar{g}_k \phi_k}{r_k^{\frac{1}{2}} \log^{\frac{\alpha}{2}} r_k } \big[z_{k+1} - G_k(\theta_k^{T}\phi_k)\big],\label{9-1}\\
    	r_k & = M^4p^2+\sum_{i=1}^{k} \bar{g}_i^2 \|\phi_i\|^2, r_0= M^4p^2,\label{9-2}
        \end{align}
\State \hspace{2em} \textbf{end for} 
\State \hspace{2em} Calculate the parameter estimates  of the mechanistic part in  (\ref{modelsimulation1}) and the bias term by using the estimate $\theta_n$:

\begin{align}
    & b_0 = \frac{\theta_n^{(2^{m_1}(1+m_2)+1)}}{\theta_n^{(1)}}, \quad  c_0 = \frac{\theta_n^{(2^{m_1}(2+2m_2)+1)}}{\theta_n^{(1)}}, \quad \bar{e} = \theta_n^{(1)}-1, \quad \quad \label{initial1*}\\
    & {p_i}_0 = \frac{\theta_n^{(i+1)}}{\theta_n^{(1)}}, \quad \quad \quad  \quad ~~{q_j}_0 = \theta_n^{(2^{m_1}j+1)}, \label{initial4*}
\end{align}
\State \textbf{\#~ Stage 2: Adam-based estimation for all unknown parameters in the model (\ref{modelsimulation1}).}

\State\hspace{2em} \textbf{for $l=1$ to $N$ do}
\State\hspace{4em} \textbf{if $l=1$  then}
\begin{equation*}\hspace{1em}
\begin{aligned}
\hat{{\Theta}}_{0 }^{(1)} &=[b_{0}, c_0, p_{10},\cdots, p_{m_10}, q_{10},\cdots, q_{m_20}, \operatorname{vec}(\hat{\Gamma}_{0})^T,  \operatorname{vec}(\hat{B}_{0})^T,  \operatorname{vec}(\hat{A}_{0})^T, (\hat{b}_{0}^{(1)})^T, (\hat{b}_{0}^{(2)})^T, \hat{b}_{0}^{(3)}]^T, 
\\ m_0^{(1)} &= 0, ~v_0^{(1)}=0.
\end{aligned}
\end{equation*}
\State\hspace{4em}  \textbf{else }

 \begin{equation*} \hspace{-21.5em}
\hat{{\Theta}}_{0}^{(l)} = \hat{{\Theta}}_{\lfloor n/\mathcal{T} \rfloor}^{(l-1)}, ~m_{0}^{(l)} = m_{\lfloor n/\mathcal{T} \rfloor}^{(l-1)},~v_{0}^{(l)} = v_{\lfloor n/\mathcal{T} \rfloor}^{(l-1)} .
\end{equation*}
\State\hspace{4em}  \textbf{end if }
\State\hspace{4em} \textbf{for $h=1$ to $\lfloor n/\mathcal{T} \rfloor$ do}
        \begin{align}
    \hat{{\Theta}}_{h}^{(l)} &=\hat{{\Theta}}_{h-1}^{(l)}
 + \frac{\eta_1}{\varepsilon + \sqrt{\frac{v_{h}^{(l)}}{1-\beta_2^{h}}}} \cdot \frac{m_{h}^{(l)}}{1-\beta_1^{h}}, \label{6}\\
     m_{h}^{(l)} &= \beta_1 m_{h-1}^{(l)} + (1-\beta_1) g_{h}^{(l)}, \label{71}\\
    v_{h}^{(l)} &= \beta_2 v_{h-1}^{(l)} + (1-\beta_2) [g_{h}^{(l)}]^2, \label{8}
\end{align}
\State \hspace{4em} \textbf{end for} 
\State \hspace{2em} \textbf{end for} 
\State  Calculate the predictive sentences in $D_2$  by using the final estimate $\hat{{\Theta}}_{\lfloor n/\mathcal{T} \rfloor}^{(N)}$. 

		\end{footnotesize}
	\end{algorithmic}
\end{algorithm*}

In Algorithm \ref{SG-Adam}, the Adam algorithm in the second stage is based on the following loss for the $h$-th batch of training data during the $l$-th epoch $(h=1,\cdots, \lfloor n/\mathcal{T} \rfloor, ~l=1,\cdots,N )$:
 \begin{equation}\label{lossadam}
    \frac{1}{\mathcal{T}} \sum_{k=(h-1)\mathcal{T}+1}^{h\mathcal{T}} \frac{|z_{k}^{(l)} - \hat{z}_{k}^{(l)}|}{z_{k}^{(l)}} + \gamma \bigg|  \frac{1}{\mathcal{T}} \sum_{k=(h-1)\mathcal{T}+1}^{h\mathcal{T}}  \hat{e}_k^{(l)}  - \bar{e}\bigg|,
\end{equation}
where  $z_k^{(l)}$ denotes the actual sentence of the $k$-th case in the $l$-th epoch, and \(\hat{z}_k^{(l)}=S_k\bigg( [ a_k + \hat{b}_{h-1}^{(l)} x_k^{(1)} + \hat{c}_{h-1}^{(l)}  x_k^{(2)} ]\)\( \times \prod_{i=1}^{m_1}\big(1 + (\hat{p}_{h-1}^{(i)})^{(l)}   v_k^{(i)}  \big) \times \big[ 1 + \sum_{j=1}^{m_2} (\hat{q}_{h-1}^{(j)} )^{(l)}  u_k^{(j)} + \hat{e}_k^{(l)} \big] \bigg)\) denotes the corresponding predictive sentence, $\hat{e}_k^{(l)} =  \hat{\Gamma}_{h-1}^{(l)} \sigma\left( \hat{B}_{h-1}^{(l)} \sigma\left( \hat{A}_{h-1}^{(l)} \eta_k + (\hat{b}_{h-1}^{(1)})^{(l)} \right) + (\hat{b}_{h-1}^{(2)})^{(l)} \right) + (\hat{b}_{h-1}^{(3)})^{(l)}$ denotes the corresponding predictive time-varying bias term, where $\hat{b}_{h-1}^{(l)}$, $\hat{c}_{h-1}^{(l)}$,  $(\hat{p}_{h-1}^{(i)})^{(l)} $,  $(\hat{q}_{h-1}^{(j)})^{(l)} $, $ \hat{\Gamma}_{h-1}^{(l)}$, $\hat{B}_{h-1}^{(l)}$, $\hat{A}_{h-1}^{(l)}$, $(\hat{b}_{h-1}^{(1)})^{(l)}$, $(\hat{b}_{h-1}^{(2)})^{(l)}$ and $(\hat{b}_{h-1}^{(3)})^{(l)}$ are the estimates of the true parameters in the model \eqref{modelsimulation1} after training on the $(h-1)$-th batch data in the $l$-th epoch.
$\bar{e}$ is the ASG-based estimate of the averaged bias term obtained by (\ref{initial1*}). The regularization coefficient $\gamma$ is determined to guide the neural network output towards matching the averaged bias term. 

Moreover, as for other related notations in Algorithm \ref{SG-Adam}, $p \triangleq 2^{m_1}(3 + 3m_2) $ in (\ref{9-2}) is the dimension of unknown parameters $\theta$ defined in (\ref{simulation3}), $\theta_n^{(\xi)} ~(\xi=1,\cdots,p)$ in (\ref{initial1*})-(\ref{initial4*}) is the $\xi$-th component of the estimate $\theta_n$. 
$g_h^{(l)}$ in (\ref{71})-(\ref{8}) denotes the negative gradient of the loss function  \eqref{lossadam} evaluated at $\hat{\Theta}_{h-1}^{(l)} =\big[\hat{b}_{h-1}^{(l)}, \hat{c}_{h-1}^{(l)}, (\hat{p}_{h-1}^{(1)})^{(l)},\cdots, (\hat{p}_{h-1}^{(m_1)})^{(l)}, (\hat{q}_{h-1}^{(1)})^{(l)}, \cdots, \\
(\hat{q}_{h-1}^{(m_2)})^{(l)}, \operatorname{vec}(\hat{\Gamma}_{h-1}^{(l)})^T,  \operatorname{vec}(\hat{B}_{h-1}^{(l)})^T,  \operatorname{vec}(\hat{A}_{h-1}^{(l)})^T, \\ \big[(\hat{b}_{h-1}^{(1)})^{(l)}\big]^T, \big[(\hat{b}_{h-1}^{(2)})^{(l)}\big]^T, (\hat{b}_{h-1}^{(3)})^{(l)}\big]^T$ with
$\operatorname{vec}(\cdot)$ denoting the vectorization of a matrix. Besides, the addition, multiplication, and division operations in the Adam algorithm are performed element-wise.

The essence of Algorithm \ref{SG-Adam} is that the ASG algorithm in the first stage generates parameter estimates close to the prediction‐error minimizer for the mechanism model, which are then used to initialize the Adam algorithm for the hybrid model in the second stage, thereby enhancing overall sentencing prediction performance.
It is worth noting that, in order to  
ensure that the neural network serves only as a compensatory component for the mechanistic model, the loss function \eqref{lossadam} of the Adam algorithm incorporates a penalty on the deviation of the averaged neural network output from the bias term $\bar{e}$ estimated by the ASG algorithm in the first stage. 




\begin{remark}
The ASG algorithm specified in (\ref{9-1})-(\ref{9-2}) is an adaptive algorithm, i.e., the algorithm updates the parameter estimate $\theta_{k+1}$ using only the current online data $\{\phi_k, z_{k+1}\}$ and the current estimate $\theta_k$. The algorithm is motivated by the classical ASG algorithm studied in \cite{2}-\cite{4}, but achieves a faster theoretical convergence rate of the averaged regret due to the use of a larger adaptation rate, which will be proven in the subsequent theoretical analysis.
\end{remark}
\begin{remark}\label{remark1}
	The estimate $\theta_k$ obtained by the ASG-type algorithm (\ref{9-1})-(\ref{9-2}) is bounded. Moreover, according to Assumption \ref{assumption3} and the definition of $\underline{g}_k, \bar{g}_k$, it follows that
	\begin{align}\label{10}
		& \inf \limits_{k \geq 0}\{\underline{g}_k\} >0,\quad  \sup \limits_{k \geq 0}\{\bar{g}_k\} <\infty.
	\end{align}
\end{remark}

\begin{remark}

For the proposed TSL algorithm, in Stage 1, the ASG method is proposed since it can provide good initial estimates for prediction in Stage 2. Specifically, the predictor based on the ASG method possesses a global asymptotic convergence property, as established in Section \ref{theory section} below, thereby providing good initial parameter estimates close to the prediction-error minimizer for the mechanism model—a guarantee typically absent in most machine learning algorithms. Moreover, the expanded regression vector in (4) is very high-dimensional in sentencing prediction, making stochastic-gradient–type methods particularly suitable for handling such cases, since the computational load of the gradient algorithms is much lower than that of the Newton-type algorithms.
In Stage 2, once the bias term $e$ is replaced by a neural network, the Adam algorithm is applied to improve the overall predictive capability of the model, as Adam is widely recognized as one of the most appropriate optimization algorithms for neural networks \cite{adamillustrate}. The advantages of using Adam are also confirmed by the experimental results presented in Section \ref{experiment section}.
\end{remark}

\section{ Prediction Theory of the ASG Algorithm} \label{theory section}


In this section, we establish the global asymptotic convergence of the ASG-based  predictor without requiring any excitation data conditions.

To facilitate the theoretical analysis, we first introduce the corresponding best predictor. By (\ref{simulation2}) and the definition of $G_k(\cdot)$, one can deduce that the best prediction of $z_{k+1}$ given $\mathscr{F}_k$ in the mean
square sense is as follows:
\begin{align}
	E[z_{k+1} | \mathscr{F}_k] = G_k(\theta^T \phi_k).
\end{align}

 Since $\theta$ is unknown \textit{a priori}, we replace the unknown parameter $\theta$ by its estimates $\theta_k$ and define the adaptive prediction for $z_{k+1}$ at time $k$ as follows:
\begin{align} \label{add3}
	\hat{	z}_{k+1}  = G_k(\theta_k^T \phi_k).
\end{align}
From the above, we define the difference between the best prediction 
and the adaptive prediction for the saturated sentence $z_{k+1}$ as “regret”, which can be expressed as follows:
\begin{align}
	R_k = [E[z_{k+1} | \mathscr{F}_k] -   \hat{	z}_{k+1}]^2 = [\psi_k]^2.
\end{align}
Naturally, we expect the regret for $z_{k+1}$ to decrease as $k$ increases, ideally vanishing to zero.
The following theorem shows that this can be achieved  in the averaged sense without requiring any data excitation condition, such as the i.i.d. condition, etc. 
\begin{theorem}\label{theorem1}
	 Under  Assumptions \ref{assumption1}-\ref{assumption3},  the accumulated  regrets have the following upper bounds:
	\begin{align}
		&	\sum_{k=0}^{n-1} R_k= o\left(n^{\frac{1}{2}} \log^{\frac{\alpha}{2}} n\right), \quad a.s.
	\end{align}
\end{theorem}
where $r_k$ and $\alpha$ are defined in Algorithm \ref{SG-Adam}. 


\begin{corollary}\label{remark2}
	If the bounded condition on $\phi_{k,i}$ in Assumption \ref{assumption1} is relaxed to $|\phi_{k,i}|  \leq  Mk^{\epsilon}, 0 < \epsilon < \frac{1}2 $ for any $1 \leq i \leq p$,
    and  $\bar{g}_k$ and $r_k$ in the algorithm (\ref{9-1})-(\ref{9-2}) are replaced by $\bar{g}_k =\sup \limits_{ {|x|\leq \max\{LMk^{\epsilon}, Mk^{\epsilon}\|\theta_k\|_1\}}}G_k^{'}(x) $ and $r_k = M^4k^{4\epsilon} p^2  +\sum_{i=1}^{k} \bar{g}_i^2 \|\phi_i\|^2  $, respectively,
	then  the accumulated regrets possess the following upper bound:
	\begin{align}
		&	\sum_{k=0}^{n-1} R_k= o\left(n^{\frac{1}{2}+\epsilon} \log ^{\frac{\alpha}{2}} n\right), \quad a.s.\label{unbound1}
	\end{align}
\end{corollary}
It can been seen that the averaged regrets $	\frac{1}{n}\sum_{k=0}^{n-1} R_k$ in both Theorem \ref{theorem1} and Corollary \ref{remark2} will converge to zero almost surely as $n \rightarrow \infty$. The proofs of the above results will be shown in the Appendix below.

\section{Sentencing Experiments} \label{experiment section}
In this section, we demonstrate the superiority of the proposed model and algorithm based on a real-world CII dataset. First, the experimental settings and the baseline models for comparison are described, followed by a detailed analysis of the results.
\subsection{Experimental Settings}
 We conduct judicial sentencing experiments based on an available CII  real-world dataset obtained from China Judgements Online\footnote{https://wenshu.court.gov.cn/}, which contains 87,588 original minor injury judgment documents, as well as 9,228 original serious injury judgment documents from the period 2019 to 2024.  These data are processed through feature extraction and quantification to facilitate subsequent experiments. 


We provide a detailed explanation of the output boundary selection and feature interpretation for the SMNN model. Regarding the output boundaries, Article 234 of Chinese Criminal Law prescribes a fixed-term imprisonment ranging from six months to three years for minor injury cases and three years to ten years for serious injury cases, so we let the statutory penalty ranges \(L_k \equiv 6\), \(N_k \equiv 36\)  and \(L_k \equiv 36\), \(N_k \equiv 120\) for minor and serious cases in (\ref{2}), respectively (Unit: Month). According to the sentencing guidelines for the CII, penalty-determining factors  $x_k^{(i)}, i=1,2$ in the model (\ref{modelsimulation1})
 represent the number of seriously and
minorly injured victims, respectively. 
Primary sentencing adjust factors $ v_k^{(i)}, i=1, \cdots,13$ in the model (\ref{modelsimulation1})  include “juveniles aged 16–18 years,” “juveniles aged 12–16 years,” “elderly individuals over 75 years of age,” “mentally disordered persons with diminished criminal responsibility,” “deaf-mute and blind individuals,” “excessive self-defense,” “excessive act of necessity,” “preparatory acts for crime,” “attempted crime,” “voluntary cessation of crime,” “accessories,” “coerced accomplices,” and “instigators.”
Other sentencing adjust factors $u_k^{(j)}, j=1, \cdots,22$ in the model (\ref{modelsimulation1}) include “grade I serious injury cases,” “grade II serious injury cases,” “grade I minor injury cases,” “grade II minor injury cases,” “voluntary surrender,” “confession,” “criminal reconciliation,” “active compensation,” “victim pardon,” “principal offender,” “crimes targeting vulnerable groups,” “victim's contributory fault,” “civil dispute-related offenses,” “recidivism,” “prior criminal records,” “first-time offenders,” “courtroom confession,” “plea agreement acceptance,” “armed affray,” “mutual combat,” “meritorious service,” and “probation,” among others.

The prediction accuracy metric adopts a relative accuracy with discretion (RAD), which is defined as follows:
\begin{equation}\label{RA}
    \text{RAD} =1- \frac{1}{n_2} \sum_{k=1}^{n_2}  \frac{\tilde{z}_k}{z_k}  I{\left(\tilde{z}_k > \max \{ 20\%z_k, 2 \} \right)},
\end{equation}
where $\tilde{z}_k =|z_k-\hat{z}_k|$ with \(z_k\) denoting the actual  sentence for the \(k\)-th judicial case, and \(\hat{z}_k\) representing its predicted value. \(n_2\) is the total number of criminal cases  in the testing set. 
The threshold  $\max \{ 20\%z_k, 2 \}$ represents the degree of adjustment allowed during the judges' sentencing process, indicating the judicial discretion.  Here the component $20\%z_k$ is based on the sentencing guidelines \footnote{See https://www.court.gov.cn/.}, while the fixed value $2$ is derived from interviews with judges and reflects practical discretion in real-world judicial decision-making.
Compared to the prediction accuracy metrics employed in previous related studies (e.g., classification metrics in \cite{41}-\cite{1}, \cite{34}), the metric (\ref{RA}) is both 
legally compatible and practically reasonable.

Moreover, the sentencing starting point \(a_k\) in the model~(\ref{modelsimulation1}) is set as the lower bound of the sentencing starting range prescribed in law based on a series of  comparative experiments(see \cite{11}). The dataset is split into a training set and a testing set with a ratio of \(4:1\). 
Some hyper-parameters for the ASG algorithm are configured as follows: the constant $\alpha$ in Algorithm~\ref{SG-Adam} is set to $1.02$,  $\mu$ and \(\bar{g}_k\) are set to \(1\), and the random noise sequence \(\{w_k\}\)  is assumed to be i.i.d. following a normal distribution \(N(0, 25)\).

Additionally, other common hyper-parameters for the Adam algorithm involved in relevant experiments are set as follows: a batch size of 245, an epoch count  of 30, a learning rate of 0.001, exponential decay rates for the first and second moment estimates  of \(\beta_1 = 0.9\), \(\beta_2 = 0.999\), a smoothing coefficient of \(\varepsilon = 10^{-8}\). The hidden layer is configured with 128 nodes.
Importantly, during training across different epochs, data are fed in chronological order without shuffling to preserve the temporal structure inherent in the dataset.
The regularization coefficient $\gamma$ in (\ref{lossadam}) is set to 0.2 and 1.4 for minor and serious injury cases respectively, based on comparative experiments of prediction accuracy.
All experiments involving random initialization are conducted 10 times, with the best-performing initialization selected to optimize overall performance.

\subsection{Baseline Models}
\subsubsection{  SM Model with ASG-Based Algorithm} \label{subsubsection1}

First, we consider the SM model (see \cite{Guo2024}), where the bias term in (\ref{modelsimulation1}) is a fixed unknown constant $e$.  Specifically, the model is as follows.
\begin{equation}\label{ex1}
	\begin{aligned}
		z_{k+1} =& S_k \bigg(\left[ a_k + b x_k^{(1)} + c x_k^{(2)} \right] \\
		&\times \prod_{i=1}^{m_1}(1 + p_i  v_k^{(i)}) \\
		&\times\left[ 1 + \sum_{j=1}^{m_2} q_j u_k^{(j)} + e \right]+ w_{k+1} \bigg).
	\end{aligned}
\end{equation}
The proposed ASG  algorithm (\ref{9-1})–(\ref{9-2}) will be used for optimization based on the model \eqref{ex1}. 

\subsubsection{SNN Model with Adam-Based Algorithm} \label{subsubsection2}
Second, we conduct experiments utilizing a saturated neural network (SNN) model (see \cite{42}) specified  as follows.
\begin{align}\label{ex2}
	&z_{k} \notag \\
    =& S_k\left[W^{(3)} \sigma\big( W^{(2)} \sigma( W^{(1)} X_k + B^{(1)} ) + B^{(2)} \big) + B^{(3)} \right],
\end{align}
where $W^{(i)}, B^{(i)}, i=1,2,3$ are unknown weighting parameters, and $X_k$ is the input data composed of various sentencing factors. Other notations are consistent with those defined in Section \ref{section2}. The Adam algorithm  (\ref{6})-(\ref{8}) will be used for  optimization based on the model \eqref{ex2}.

\subsubsection{ SMNN Model with Adam-Based Algorithm} \label{subsubsection3}
Third,  we train  our SMNN model  (\ref{modelsimulation1}) by Adam optimizer  (\ref{6})–(\ref{8}) with all the parameters being initialized randomly.    

\subsubsection{ SMNN Model with TSL Algorithm} 
\label{subsubsection4}
Fourth,  we train the SMNN model  (\ref{modelsimulation1}) by our TSL optimizer described in Algorithm \ref{SG-Adam}, where the mechanistic parameters are initialized using the ASG algorithm in the first stage, and other parameters are initialized randomly.

\subsection{Comparison Experiments on Sentencing Prediction}
To further improve the sentencing prediction accuracy for CII in the existing literature (e.g., \cite{11},  \cite{Guo2024}, \cite{newmodel}), we perform comparison experiments to demonstrate the improved accuracy of our SMNN model and the TSL algorithm compared to the above baseline methods using the same available data. 

\begin{figure}
		\includegraphics[width=9cm,height=6cm]{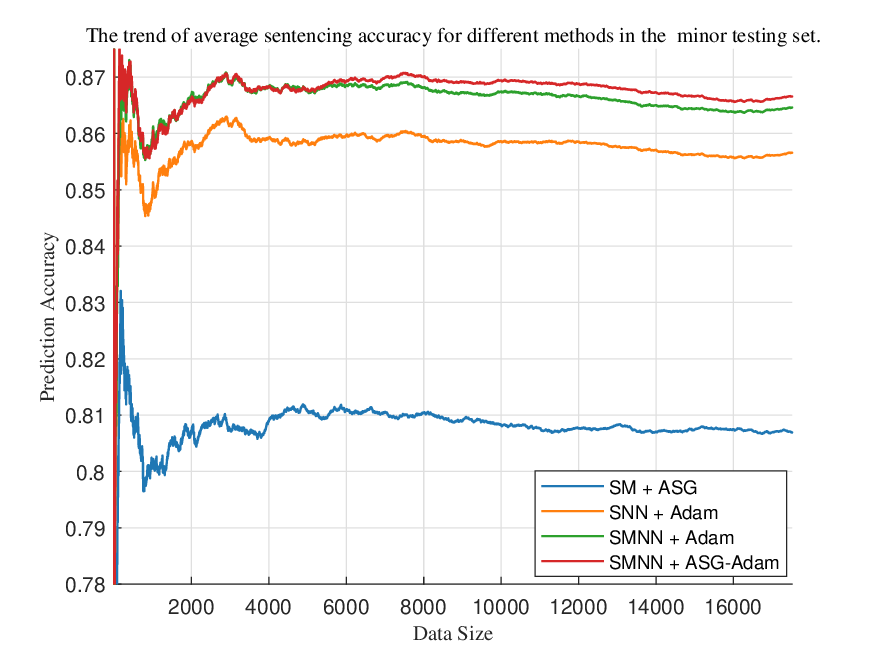}    
		\caption{ The trend of average sentencing accuracy of minor cases for different methods in the testing set.}  
		\label{fig1}                                 
\end{figure}

\begin{figure}
		\includegraphics[width=9cm,height=6cm]{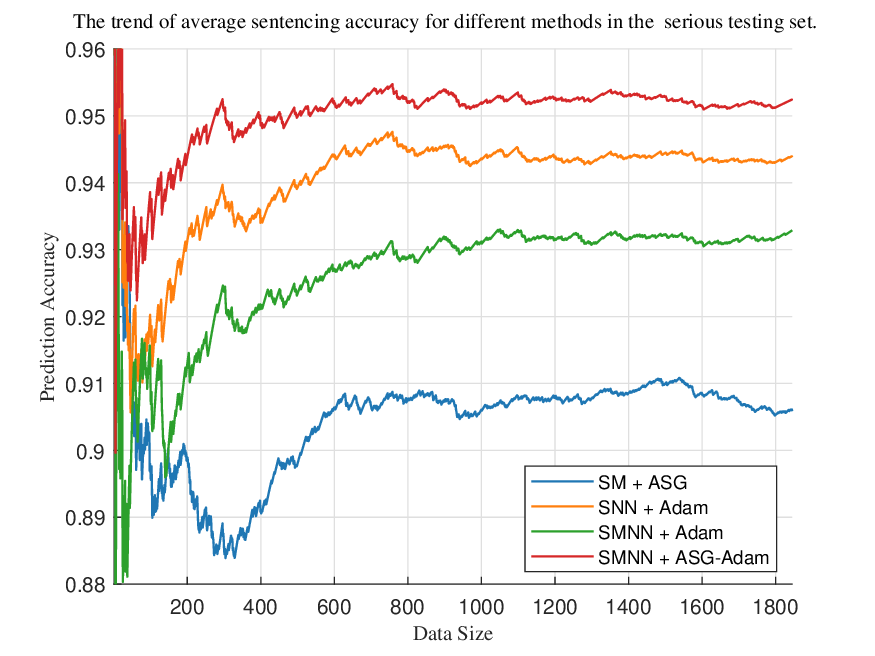}    
		\caption{ The trend of average sentencing accuracy of serious cases for different methods in the testing set.} 
		\label{fig2}                                 
\end{figure}

Figs. \ref{fig1}-\ref{fig2}  illustrate the prediction accuracy trends of  our method in comparison to other related known methods on the testing set for minor and serious cases, respectively. It can be shown that:
\begin{itemize}
    \item The SMNN model with the ASG-Adam algorithm consistently achieves the highest prediction accuracy of 86.66\% and 95.25\% for minor and serious cases respectively, surpassing the SM model with the ASG algorithm (80.69\%; 90.59\%) and the SNN model with the Adam algorithm (85.66\%; 94.38\%). This demonstrates the advantages of combining both the SM model with the NN model as well as integrating the ASG algorithm with the Adam algorithm, allowing for more accurate sentencing prediction.
    \item The SMNN model with  our ASG-Adam algorithm outperforms the same model with the Adam algorithm with  random initialization (86.46\%; 93.29\%), which demonstrates that our ASG algorithm (\ref{9-1})-(\ref{9-2}) in the first stage provides good initial parameter estimates for the Adam optimizer (\ref{6})-(\ref{8}) in the second stage, thereby improving prediction accuracy. This advantage stems from the ASG-based initialization, which guides the optimization process toward the global error minimum. In contrast, random initialization requires exploring a broader solution space and is more prone to local error minima, often resulting in inferior prediction accuracy.    
    
\end{itemize}


\section{Conclusion} 
Motivated by the imperative demand in judicial practice for highly accurate and reliable sentencing prediction, this paper proposed a SMNN hybrid model, which integrates the sentencing logic model  with a neural network. 
The ASG algorithm is applied to update the parameters of the  mechanistic component within the SMNN model, after which the Adam algorithm is utilized to further optimize the SMNN model with the aim of enhancing predictive accuracy.
Both theoretical analysis and sentencing experiments demonstrate that the estimates produced by the ASG algorithm in the first stage can provide good enough initial values for the implementation of the prediction algorithm in the second stage.
Moreover, sentencing experiments also reveal that the neural networks introduced in the hybrid model is indeed helpful for improving the prediction accuracy in the second stage.
For future investigation, it would  be interesting to establish the global convergence of the Adam algorithm in the second stage, and to apply our prediction algorithms to other judicial crimes other than the intentional injury studied here.


\section{Appendix}
\begin{lemma}\label{lemma1}  {\rm (\cite{6} Theorem 1.3.2)}\ \
	Let $\{f_k,\mathscr{F}_k\}$ and $\{\alpha_k,\mathscr{F}_k\}$ be two non-negative adapted sequences. If $E[f_{k+1}| \mathscr{F}_k] \leq f_k +\alpha_k,  a.s.$ for any $k \geq 1$ and $\sum_{i=1}^{\infty} \alpha_i < \infty, a.s.$, then $f_k$ will converge to a  finite limit a.s.
\end{lemma}
\begin{lemma} \label{lemma2} {\rm (\cite{6} Theorem 1.2.15)} \ \
	Let $D_k = C + \sum_{j=1}^{k} d_j, d_j \geq 0$ and $D_0=C$ with $C>1$ being any constant, then we have
	\begin{align}\label{3.1}
		& \sum_{j=1}^{\infty} \frac{d_j}{D_j \log^{\alpha} D_j } < \infty, \forall \alpha >1.
	\end{align}
\end{lemma}

\begin{lemma}  \label{lemma8} {\rm (\cite{9})}\ \
	Let $\left\{\omega_n, \mathcal{F}_n\right\}$ be a martingale difference sequence and $\left\{f_n, \mathcal{F}_n\right\}$ an adapted sequence. If
	\begin{align}
		\sup _n \mathbb{E}\left[\left|\omega_{n+1}\right|^\alpha \mid \mathcal{F}_n\right]<\infty, \text { a.s. }
	\end{align}
	for some $\alpha \in(0,2]$, then as $n \rightarrow \infty$, we  have $\forall \eta>0$,
	\begin{align}
		\sum_{i=0}^n f_i \omega_{i+1}=O\left(s_n(\alpha) \log ^{\frac{1}{\alpha}+\eta}\left(s_n^\alpha(\alpha)+e\right)\right) \text { a.s., }
	\end{align}
	where
	$
	s_n(\alpha)=\left(\sum \limits_{i=0}^n\left|f_i\right|^\alpha\right)^{\frac{1}{\alpha}}.
	$
\end{lemma}
Inspired by the analysis ideas mentioned in \cite{6},\cite{12}, and \cite{21}, etc, the proof of Theorem \ref{theorem1} is as follows:

\textbf{Proof of Theorem \ref{theorem1}. } 
By the algorithm  (\ref{9-1})-(\ref{9-2}) and the differential mean value theorem, we have
\begin{align}
	\tilde{\theta}_{k+1}&=\tilde{\theta}_k-\frac{\mu \bar{g}_k \phi_k}{r_k^{\frac{1}{2}} \log^{\frac{\alpha}{2}} r_k}\left[G_k^{'}(\xi_k)\phi_k^{T} \tilde{\theta}_k+v_{k+1}\right]  \notag\\
	& =\left(I - \frac{\mu \bar{g}_k G_k^{'}(\xi_k)  \phi_k  \phi_k^{T}}{r_k^{\frac{1}{2}} \log^{\frac{\alpha}{2}} r_k}\right)\tilde{\theta}_k - \frac{\mu \bar{g}_k \phi_k}{r_k^{\frac{1}{2}} \log^{\frac{\alpha}{2}} r_k} v_{k+1}, \label{2.2}
\end{align}
where $\xi_k  \in (\min\{\theta^\tau \phi_k, \theta_k^\tau \phi_k\},\max\{\theta^\tau \phi_k, \theta_k^\tau \phi_k\})$.
By this and the elementary inequality we have
\begin{align}
	&\tilde{\theta}_{k+1}^T \tilde{\theta}_{k+1} \notag \\
	= &\left\{\left(I - \frac{\mu \bar{g}_k G_k^{'}(\xi_k)  \phi_k  \phi_k^{T}}{r_k^{\frac{1}{2}} \log^{\frac{\alpha}{2}} r_k}\right)\tilde{\theta}_k - \frac{\mu \bar{g}_k \phi_k}{r_k^{\frac{1}{2}} \log^{\frac{\alpha}{2}} r_k} v_{k+1}\right\}^T \notag \\
	&\left\{\left(I - \frac{\mu \bar{g}_k G_k^{'}(\xi_k)  \phi_k  \phi_k^{T}}{r_k^{\frac{1}{2}} \log^{\frac{\alpha}{2}} r_k}\right)\tilde{\theta}_k - \frac{\mu \bar{g}_k \phi_k}{r_k^{\frac{1}{2}} \log^{\frac{\alpha}{2}} r_k} v_{k+1}\right\} \notag \\
	= & \tilde{\theta}_{k}^T \tilde{\theta}_{k} - 2 \frac{\mu \bar{g}_k G_k^{'}(\xi_k) \left[\tilde{\theta}_k^T  \phi_k \right]^2}{r_k^{\frac{1}{2}} \log^{\frac{\alpha}{2}} r_k } - 2 \frac{\mu \bar{g}_k \tilde{\theta}_k^T  \phi_k}{r_k^{\frac{1}{2}} \log^{\frac{\alpha}{2}} r_k }v_{k+1} \notag \\
	&+\frac{\mu^2\bar{g}_k^2 [G_k^{'}(\xi_k)]^2 \|\phi_k\|^2 \left[\tilde{\theta}_k^T  \phi_k \right]^2 }{r_k \log ^{\alpha}r_k } \notag \\
	& + 2 \frac{ \mu^2 \bar{g}_k^2 G_k^{'}(\xi_k) \|\phi_k\|^2 \tilde{\theta}_k^T  \phi_k  }{r_k \log ^{\alpha}r_k } v_{k+1} + \frac{ \mu^2 \bar{g}_k^2 \|\phi_k\|^2}{r_k \log ^{\alpha}r_k}  v_{k+1}^2 \notag \\
	\leq & \tilde{\theta}_{k}^T \tilde{\theta}_{k} -2 \frac{\mu \bar{g}_k G_k^{'}(\xi_k)\left[\tilde{\theta}_k^T  \phi_k \right]^2}{r_k^{\frac{1}{2}} \log^{\frac{\alpha}{2}} r_k } -2 \frac{\mu  \bar{g}_k \tilde{\theta}_k^T  \phi_k}{r_k^{\frac{1}{2}} \log^{\frac{\alpha}{2}} r_k }  v_{k+1} \notag \\
	&+ (1+b) \frac{ \mu^2 \bar{g}_k^2 [G_k^{'}(\xi_k)]^2 \|\phi_k\|^2 \left[\tilde{\theta}_k^T  \phi_k \right]^2 }{r_k \log^{\alpha} r_k} \notag\\
	&+ \left(1+ \frac{1}{b}\right)\frac{ \mu^2 \bar{g}_k^2 \|\phi_k\|^2}{r_k \log^{\alpha} r_k} v_{k+1}^2 \label{p1},
\end{align}
where $0<b<1$ is a constant. since $S_k^{'}(x) \leq 1, \forall x \in \mathbb{R}$ implies $0 \leq  \bar{g}_k G_k^{'}(\xi_k) \leq 1$, summing up both sides of (\ref{p1}) from $k=0$ to $n-1$, and noticing $\|\phi_k\|^2 < r_k^{\frac{1}{2}} \log ^{\frac{\alpha}{2}} r_k$  for all $k \geq 0$, we know that
\begin{align} \label{p2}
	&\| \tilde{\theta}_{n}  \|^2 + (1-b) \sum_{k=0}^{n-1}  \frac{ \mu  \bar{g}_k G_k^{'}(\xi_k)  \left[\tilde{\theta}_k^T  \phi_k \right]^2}{ r_k^{\frac{1}{2}} \log ^{\frac{\alpha}{2}} r_k } \notag \\
	= & O\left( \sum_{k=0}^{n-1} \frac{ \bar{g}_k\tilde{\theta}_k^T  \phi_k}{r_k^{\frac{1}{2}} \log^{\frac{\alpha}{2}} r_k }  v_{k+1}  \right) + O\left( \sum_{k=0}^{n-1}  \frac{ \bar{g}_k^2\|\phi_k\|^2}{r_k \log^{\alpha} r_k} v_{k+1}^2 \right).
\end{align}
Now we analyze the RHS of (\ref{p2}) term by term. Note that  by (\ref{7}), the elementary inequality $ 2ab \leq a^2 + b^2$ and the inequality $ E^{4}[|x||\mathscr{F}_{k}] \leq E[|x|^4|\mathscr{F}_{k}]$, we have
\begin{align}\label{add8} 
	& E\big[\left|v_{k+1}\right|^{4} \mid \mathscr{F}_{k}\big] \notag \\
	=& O\bigg(E\big[\big|S_k(\theta^{\tau} \phi_k +w_{k+1})-S_k(\theta^{\tau} \phi_k)\big|^4 \big| \mathscr{F}_{k} \big]\bigg) \notag \\
	 + &O\bigg(E\big[\big|S_k(\theta^{\tau} \phi_k)-E[ S_{k}(\theta ^{\tau} \phi_k +w_{k+1}) \big|\mathscr{F}_{k}] \big|^4 \big| \mathscr{F}_{k} \big]\bigg) \notag \\
	=& O\bigg(E\big[\big|S_k(\theta^{\tau} \phi_k +w_{k+1})-S_k(\theta^{\tau} \phi_k)\big|^4 \big| \mathscr{F}_{k} \big]\bigg)+O(1) \notag \\
	=& O\bigg(E\big[|w_{k+1}|^4 \big|\mathscr{F}_{k}\big]\bigg)+O(1)= O(1) , \quad a.s.
\end{align}
Note that by the definition of $v_{k+1}$, we have
\begin{align}\label{1.4}
	& E\left(v_{k+1} \mid \mathscr{F}_k\right) \notag \\
	=&  E\left[S_k\left(\theta^\tau \phi_k+w_{k+1}\right) \mid \mathscr{F}_k\right] \notag \\
	  & -E\left[E\left[S_k\left(\theta^\tau \phi_k+w_{k+1}\right) \mid \mathscr{F}_k\right] \mid \mathscr{F}_k\right] =0.
\end{align}
So by (\ref{add8}),  (\ref{1.4}) and Lemma \ref{lemma8}, we have
\begin{align}\label{p4}
	\sum_{k=0}^{n-1} \frac{ \bar{g}_k\tilde{\theta}_k^T  \phi_k}{r_k^{\frac{1}{2}} \log^{\frac{\alpha}{2}} r_k }  v_{k+1} & = o\left( \sum_{k=0}^{n-1}  \frac{  \bar{g}_k^2 \left[\tilde{\theta}_k^T  \phi_k \right]^2}{r_k^{\frac{1}{2}} \log^{\frac{\alpha}{2}} r_k}\right)+ O(1),  a.s.
\end{align}
Moreover, by  (\ref{add8}), (\ref{1.4}) and Lemma \ref{lemma2} and Lemma \ref{lemma8}, we have
\begin{align}\label{p5}
	& \sum_{k=0}^{n-1}  \frac{ \bar{g}_k^2\|\phi_k\|^2}{r_k \log^{\alpha} r_k} v_{k+1}^2  \notag \\
    \leq &	\sum_{k=0}^{n-1}  \frac{ \bar{g}_k^2 \|\phi_k\|^2}{r_k \log^{\alpha} r_k} \left[v_{k+1}^2 - E[v_{k+1}^2|\mathscr{F}_k]\right] \notag\\
	& + 	 \sup_{k \geq 0} E[v_{k+1}^2|\mathscr{F}_k] \cdot\sum_{k=0}^{n-1}  \frac{ \bar{g}_k^2 \|\phi_k\|^2}{r_k \log^{\alpha} r_k} \notag \\
	 = &  O\left(\sum_{k=0}^{n-1}  \frac{  \bar{g}_k^2 \|\phi_k\|^2}{r_k \log^{\alpha} r_k}\right) + O(1) \notag\\
	 = & O(1), \quad  a.s.
\end{align}
Combining (\ref{p4})-(\ref{p5}) with (\ref{p2}), we have
\begin{align} \label{p6}
	\sum_{k=0}^{n-1}  \frac{  \bar{g}_k G_k^{'}(\xi_k)  \left[\tilde{\theta}_k^T  \phi_k \right]^2}{ r_k^{\frac{1}{2}} \log ^{\frac{\alpha}{2}} r_k }  = O(1), \quad a.s.
\end{align}
Note that since $\phi_k$ is bounded, we have $r_k = O(k)$. 
Moreover, if $r_k \rightarrow \infty$ as $ k \rightarrow \infty$, then we have
by (\ref{p6})  and the Kronecker Lemma, we obtain
\begin{align}\label{*11}
	\sum_{k=0}^{n-1} R_k  = o\left(n^{\frac{1}{2}} \log^{\frac{\alpha}{2}} n\right), \quad a.s.
\end{align}
Otherwise, if $r_k \nrightarrow \infty$, then  denominator of (\ref{p6}) is bounded, and (\ref{*11}) is also satisfied.

\textbf{Proof of Remark \ref{remark1}.} 
By (\ref{9-1}), we have
\begin{align}\label{1.1}
	\tilde{\theta}_{k+1}^{T} \tilde{\theta}_{k+ 1}=&\left[\tilde{\theta}_k-\frac{ \mu \bar{g}_k \phi_k}{r_k^{\frac{1}{2}} \log ^{\frac{\alpha}{2}} r_k}\left[z_{k+1}-G_k\left(\phi_k^{T} \theta_k\right)\right]\right]^{T} \notag\\
	& \left[\tilde{\theta}_k-\frac{\mu \bar{g}_k \phi_k}{r_k^{\frac{1}{2}} \log ^{\frac{\alpha}{2}} r_k}\left[z_{k+1}-G_k\left(\phi_k^{T} \theta_k\right)\right]\right].
\end{align}
Set $\sup \limits _{k \geq 0} \max\{|L_k|,|N_k|\} \leq  U < \infty$, from (\ref{26}) and (\ref{27}) we know that
\begin{align}
	& \left\|\tilde{\theta}_{k+1}\right\|^2 \notag \\
    = & \left\|\tilde{\theta}_k\right\|^2-2 \frac{\mu \bar{g}_k}{r_k^{\frac{1}{2}} \log ^{\frac{\alpha}{2}} r_k}\left[z_{k+1}-G_k\left(\phi_k^{T} \theta_k\right)\right] \tilde{\theta}_k^{T} \phi_k \notag  \\
	&+\frac{\mu^2\bar{g}_k^2\left[z_{k+1}-G_k\left(\phi_k^{T}\theta_k\right)\right]^2}{r_k \log ^{\alpha} r_k}\left\|\phi_k\right\|^2 \notag \\
	 \leq & \left\|\tilde{\theta}_k\right\|^2-2\frac{\mu\bar{g}_k}{r_k^{\frac{1}{2}} \log ^{\frac{\alpha}{2}} r_k}\left[v_{k+1}+\psi_k\right] \tilde{\theta}_k^{T} \phi_k \notag\\
	& +\frac{4U^2\bar{g}_k^2}{r_k \log^{\alpha} r_k}\left\|\phi_k\right\|^2. \label{1.3}
\end{align}
Take the conditional expectation for both sides of (\ref{1.3}), and by (\ref{1.4}) and differential mean value theorem, we know that there exists a random variable $\xi_k \in (\min\{\theta^\tau \phi_k, \theta_k^\tau \phi_k\},\max\{\theta^\tau \phi_k, \theta_k^\tau \phi_k\})$ such that
\begin{align}
	& E\left[\left\|\tilde{\theta}_{k+1}\right\|^2 | \mathscr{F}_k\right]  \notag \\
    \leq & \left\|\tilde{\theta}_k\right\|^2-2\frac{\mu \bar{g}_kG_k^{'}(\xi_k) (\tilde{\theta}_k^{T} \phi_k)^2 }{r_k^{\frac{1}{2}} \log ^{\frac{\alpha}{2}} r_k} +\frac{4U^2\bar{g}_k^2\left\|\phi_k\right\|^2 }{r_k \log^{\alpha} r_k}\notag\\
	 \leq & \left\|\tilde{\theta}_k\right\|^2+\frac{4U^2\bar{g}_k^2\left\|\phi_k\right\|^2}{r_k \log^{\alpha} r_k}. \label{1.6}
\end{align}
Moreover, by Lemma \ref{lemma2} we obtain $\sum \limits_{k=1}^{n} \frac{\bar{g}_k^2\|\phi_k\|^2}{r_k \log^{\alpha} r_k} = O(1)$, by this, (\ref{1.6}) and  Lemma \ref{lemma1},  we know that there exists a constant $S  \geq  0$ such that  $\|\tilde{\theta}_k\| \rightarrow S < \infty $ as $k \rightarrow \infty$. Combining this with the boundedness of $\theta$, we can easily obtain the estimate $\theta_k$ is bounded.

\textbf{Proof of Corollary  \ref{remark2}.}  Note that  $\|\phi_k\|^2 \leq r_k^{\frac{1}{2}} \log ^{\frac{\alpha}{2}} r_k $, similar to the analysis of (\ref{p2}) and (\ref{p4})-(\ref{p5}), we have
\begin{align}
	\sum_{k=0}^{n-1} \bar{g}_k G_k^{'}(\xi_k)  \left[\tilde{\theta}_k^T  \phi_k \right]^2 = o(n^{\frac{1}{2}+\epsilon} \log ^{\frac{\alpha}{2}} n), a.s.,
\end{align}
which implies that (\ref{unbound1}) holds.

\newpage

\vfill

\end{document}